\documentclass[12pt]{amsart}
\title {Combinatorics of the heat trace on spheres} 
\author{Iosif Polterovich}
\thanks{Supported by CRM-ISM and MSRI 
postdoctoral fellowships}
\address{Centre de Recherches Math\'ematiques, Universit\'e de Montr\'eal,
C.P. 6128 Succ. centre-ville, Montr\'eal (Qu\'ebec) H3C 3J7, Canada}
\email{iossif@math.uqam.ca}  
\def \phi{\varphi}
\def \epsilon{\varepsilon}
\numberwithin{equation}{subsection}
\theoremstyle{definition}

\theoremstyle{plain}
\newtheorem{lemma}[equation]{Lemma}
\newtheorem{theorem}[equation]{Theorem}

\begin{document}
\maketitle
%\centerline{\large \bf Preliminary version}
\begin{abstract} 
We present a concise explicit expression for the heat trace coefficients
of spheres. Our formulas yield certain combinatorial
identities which are proved following ideas of D.~Zeilberger. In particular,
these identities  allow to recover in a surprising way some known formulas 
for the heat  trace asymptotics.  
Our approach is based on a method for computation of heat invariants 
developed in [P].
\end{abstract}
\section{Introduction and main results}
\subsection{Heat trace asymptotics on spheres}
Let $S^d$ be a sphere with the standard Riemannian metric of curvature $+1$.
The Laplace-Beltrami operator $\Delta$ on $S^d$ has eigenvalues 
$\lambda_{k,d}=k(k+d-1)$, and each $\lambda_{k,d}$ has multiplicity
$\mu_{k,d}$ given by
$$
\mu_{k,d}=\frac{(2k+d-1)(k+d-2)!}{k! (d-1)!},\, \, k\ge 1 \,\,\,
\operatorname{and} \,\, \mu_{0,d}=1,
$$ 
(see [M\"u]).
Consider an asymptotic expansion for the trace of the heat operator 
$e^{-t\Delta}$ as $t\to 0+$ (see [Be], [Gi]):
\begin{equation}
\label{tr}
\sum_{\lambda} e^{-t\lambda}=\sum_{k=0}^{\infty} 
\mu_{k,d}\, e^{-t\lambda_{k,d}} 
\sim  \sum_{n=0}^{\infty}a_{n,d} t^{n-\frac{d}{2}}.
\end{equation} 
Heat trace coefficients (or heat invariants) $a_{n,d}$ were  
calculated in [CW] (see (\ref{odd}) and (\ref{even}) for similar formulas)
by methods of Lie groups and representation theory 
(see also [Ca], [ELV], [DK] for related results). 
In this paper we present a different approach based on [P]. We obtain the 
following concise explicit expression for $a_{n,d}$.
\begin{theorem}
\label{s}
For any $n\ge 1$ and any integer $\omega \ge 2n$ the heat invariants 
$a_{n,d}$ are equal to
\begin{multline}
\label{sp}
a_{n,d}=
\sum_{j=1}^{\omega}
\frac{2 (-1)^n \Gamma(\omega +\frac{d}{2}+1)}
{(\omega-j)!(j+n)!(2j+d)!}
\sum_{k=1}^j (-1)^k \binom{2j+d-1}{j-k}\mu_{k,d}\, \lambda_{k,d}^{j+n} 
\end{multline}
\end{theorem}
There is some delicacy in the proof of Theorem \ref{s}. For $\omega \ge 3n$
it follows from a simple generalization of the main result of [P] and
some facts about Legendre polynomials (see sections 2.1 and 2.2).
Theorem \ref{s} for $2n\le \omega <3n$ follows from the proofs of Theorems
\ref{o} and \ref{e} involving rather sophisticated combinatorial 
arguments  due to Doron Zeilberger (see below). 

Validity of formula (\ref{sp}) for $3n>\omega \ge 2n$ was suggested by
computer experiments using [Wo]. Note that $2n$ is ``sharp''
in a sense that if $\omega < 2n$ then (\ref{sp}) is no longer true 
(see section 3.1).
\subsection{Combinatorial identities}
Taking $d=1$ in (\ref{sp}) we should get zero since the heat trace
coefficients $a_{n,1}$ of a circle $S^1$ vanish identically for $n\ge 1$.
This gives rise to a surprising combinatorial identity:
\begin{theorem} ($S^1$-identity) \mbox{\rm (D.~Zeilberger, [Z])}
\label{S1}
$$
\sum_{j=0}^{\omega} \frac{1}{(\omega-j)! (j+n)! (2j+1)}\sum_{k=0}^j
\frac{(-1)^k k^{2j+2n}}{(j-k)!(j+k)!} = 0
$$
for $n\ge 1$, $\omega\ge 2n$.
\end{theorem}
Theorem \ref{S1} was proved in [Z] (see also section 3.1)
by pure combinatorial methods.

Similarly, taking into account that
$$a_{n,3}=\frac{\sqrt{\pi}}{4\cdot n!}\,\,\operatorname{(cf. [MS], [CW])},
$$
we get
\begin{theorem} ($S^3$-identity)
\label{S3}
$$
\sum_{j=0}^{\omega} \frac{\Gamma(\omega+5/2)}
{(\omega-j)! (j+n)! (2j+3)} \sum_{l=0}^{j+1} 
\frac{(-1)^{l} l^2 (l^2-1)^{j+n}}{(j+l+1)!(j-l+1)!}=
\frac{(-1)^{n+1}\sqrt{\pi}}{8\cdot n!},
$$
for $n\ge 1$, $\omega\ge 2n$.
\end{theorem}
A  combinatorial  proof of this theorem based on a generalization of 
Zeilberger's arguments is given in section 3.2.

Interestingly enough, pushing forward this combinatorial approach 
one recovers the results of [CW] from Theorem \ref{s} for $\omega \ge 3n$.
We present them in a more concise form especially in some particular cases 
(see (\ref{57}), (\ref{S2})).
\subsection{Odd-dimensional case}
In odd dimensions formula (\ref{sp}) can be substantially simplified.
\begin{theorem}
\label{o}
The heat invariants of odd-dimensional spheres $S^{2\alpha +1}$
are equal to
\begin{equation}
\label{odd}
a_{n,2\alpha+1}=\sum_{s=1}^{\alpha} 
\frac{\alpha^{2n-2\alpha+2s}\Gamma(s+\frac{1}{2})
K_s^{\alpha}}{(n-\alpha+s)! (2\alpha)!},
\end{equation}
where the coefficients $K_s^{\alpha}$ are defined by
\begin{equation}
\label{cno}
\prod_{\beta=0}^{\alpha-1} (z^2-\beta^2) = \sum_{s=1}^{\alpha} 
K_s^{\alpha} z^{2s}.
\end{equation}
In particular, 
\begin{equation}
\label{57}
a_{n,5}=\frac{4^{n-3} (6-n)\sqrt{\pi}}{3\cdot n!},\quad
a_{n,7}=\frac{3^{2n-6} (16n^2-286n+1215)\sqrt{\pi}}{640\cdot n!}.
\end{equation}
\end{theorem}
\subsection{Even-dimensional case}
Formulas for even-dimensional spheres have a more intricate
combinatorial structure  due to a certain hypergeometric expression 
vanishing only for $d$ odd (see section 4.2).
\begin{theorem}
\label{e}
The heat invariants of even-dimensional spheres $S^{2\nu}$
are equal to
\begin{multline}
\label{even}
a_{n,2\nu}=
\frac{1}{(2\nu-1)!}\left(\sum_{t=0}^{\nu-1} \frac{(\nu-1-t)!}{(n-t)!}
\left(\nu-\frac{1}{2}\right)^{2n-2t}K_t^\nu +\right.\\
\left.\sum_{t=0}^{\nu-1}K_t^\nu \sum_{p=\nu-t}^{n-t} (-1)^{p+\nu-t-1}
\frac{(\nu-\frac{1}{2})^{2n-2t-2p} B_{2p}}{p\, (n-t-p)! (p-\nu+t)!}
\left(\frac{1}{2^{2p-1}}-1\right)\right)
\end{multline}
where $B_{2p}$ are the  Bernoulli numbers (see {\rm [GKP]}) and the 
constants $K_t^\nu$
are defined by
\begin{equation}
\label{cne}
\prod_{\beta=1/2}^{\nu-3/2} (z^2-\beta^2) = \sum_{t=0}^{\nu-1}K_t^\nu z^{2\nu-2-2t}
\end{equation}
\end{theorem}
In particular (cf. [Ca]),
\begin{equation}
\label{S2}
a_{n,2}=\frac{1}{n! 2^{2n}}\sum_{r=0}^n(-1)^r\binom{n}{r}
(2-2^{2r})B_{2r}.
\end{equation}
Note that the second sum in (\ref{even}) vanishes for $\nu>n$.
\subsection{Structure of the paper}
In section 2.1 we present a generalization of the main result
of [P] which allows to prove Theorem \ref{s} for $\omega\ge 3n$ 
using some properties of Legendre polynomials, see section 2.2. 
In section 3.1 we review Zeilberger's proof of  Theorem \ref{S1} 
which leads to the  proof of Theorem \ref{S3} in section 3.2. 
Theorems \ref{o} and \ref{e} are proved in  sections 4.1 and 4.2
using Theorem \ref{s} for $\omega \ge 3n$.
Theorem \ref{s} for $3n > \omega \ge 2n$ 
follows from Theorems \ref{o} and \ref{e} by reversing arguments in 
their proofs, see section 4.3. 
Two auxiliary combinatorial lemmas are proved in sections 5.1 and 5.2.
\subsection*{Acknowledgments} I am very grateful to Doron Zeilberger for
helpful advice  concerning combinatorial identities proved
in this paper, and especially for his proof of Theorem \ref{S1}.
I would like to thank Dmitry Jakobson, Yakar Kannai, Leonid Polterovich, 
Andrei Reznikov and Joseph Wolf for stimulating  discussions.
The author is also indebted to Klaus Kirsten and Francois Lalonde 
for useful remarks on the first draft of this paper.

This research was partially conducted during my stay at the Mathematical
Sciences Research Institute in Berkeley whose hospitality and support
are gratefully acknowledged.

\section{Heat invariants and spherical harmonics}
\subsection{Computation of heat invariants}
For any $d$-dimensional closed Riemannian manifold $M$ the coefficients 
$a_{n,d}$ can be obtained from the local heat invariants $a_{n,d}(x)$ 
(see [B], [Gi], [P]): 
$$a_{n,d}=\int\limits_M a_{n,d}(x) d\operatorname{vol}(x).$$
In particular, if $M=S^d$ the coefficients $a_{n,d}(x)$ are constants and 
therefore for any $x\in S^d$
\begin{equation}
\label{vol1}
a_{n,d}=\operatorname{vol}(S^d)a_{n,d}(x),
\end{equation}
where the volume of a $d$-sphere is given by (see [M\"u]):
\begin{equation}
\label{vol2}
\operatorname{vol}(S^d)=\frac{2\pi^{\frac{d+1}{2}}}{\Gamma(\frac{d+1}{2})}.
\end{equation}
Let us prove the following modification of the main result of [P]:
\begin{theorem}
\label{most1}
For any integer $\omega \ge 3n$ the local heat invariants $a_{n,d}(x)$ 
of a $d$-dimensional closed Riemannian manifold $M$ are equal to: 
\begin{equation}
\label{mm}
a_{n,d}(x)=(4\pi)^{-d/2}(-1)^n\sum_{j=0}^{\omega}
\binom{\omega+\frac{d}{2}}{j+\frac{d}{2}}
\frac{1}{4^j \, j! \, (j+n)!}\left.
\Delta^{j+n}(f(r_x(y)^2)^j)\right|_{y=x},
\end{equation}
where $f(r_x^2)$ is a smooth function in some 
neighborhood of $x\in M$ such that $f(s)=s+O(s^2)$, $s\in [0,\epsilon]$.
\end{theorem}
\noindent {\bf Proof.}
The result follows from Theorem 1.2.1 in [P]) (if $f(r_x^2)=r_x^2$ 
we get precisely the statement of that theorem). Indeed,
let $(y_1,\dots,y_d)$ be normal coordinates in a neighborhood of 
the point $x=(0,\dots,0)\in M$. The Riemannian metric at the point $x$ 
has the form $ds^2=dy_1^2+\cdots + dy_d^2$ and 
the square of the distance function is locally given by
\begin{equation}
\label{1}
r_x(y)^2=y_1^2+\cdots+y_d^2,
\end{equation}
where $y=(y_1,\dots,y_d)$. 
Let us note that the point $x\in M$ is a 
non-degenerate critical point of index $0$ of
the function $f$ and hence due to Morse lemma ([Mi])
the function $f(r_x^2(y))$ can be locally written as the sum of squares
(\ref{1}) in some new coordinate system $(y_1',\dots,y_d')$. 
Moreover, this new system can be chosen in such a 
way that $y_1'=y_1+O(|y|^2),\dots,y_d'=y_d+O(|y|^2)$, and hence
the Riemannian metric remains Euclidean at the point $x$.
Repeating the proof  of Theorem 1.2.1 in [P] with the 
coordinates  $(y_1',\dots, y_d')$ taken instead of normal coordinates 
we complete the proof of (\ref{mm}). \qed

\noindent {\bf Remark.} As was recently observed in [We], for
$f(r^2)=r^2$ one could in fact take $\omega \ge n$ in (\ref{mm}).
This can be also deduced from the proof of Theorem 1.2.1 in [P]
taking into account  that the radial part of the Laplacian in  normal 
coordinates is a first order perturbation of the operator 
$-\partial^2/\partial r^2$.
\subsection{Application of Legendre polynomials}
Recall that the Laplacian on $S^d$ has eigenvalues 
$\lambda_{k,d}=k(k+d-1)$ and 
the corresponding eigenfunctions are the Legendre polynomials 
$L_{k,d}(\cos r)$ (see [M\"u]):
\begin{equation}
\label{leg}
\Delta L_{k,d}=\lambda_{k,d} L_{k,d}=k(k+d-1) L_{k,d}
\end{equation} 
\noindent{\bf Proof of Theorem \ref{s} for $\omega \ge 3n$.}
Take  $f(r^2)=2-2\cos(r)=r^2+O(r^4)$ as the function $f$ in 
Theorem \ref{most1}. We express its powers in terms of the Legendre 
polynomials $L_{k,d}(\cos r)$. Denote $t=\cos r$. Then $f(r^2)^j=2^j(1-t)^j$.
Let 
\begin{equation}
\label{leg0}
f(r^2)=2^j (1-t)^j=2^j\sum_{k=0}^j c_{jk}L_{k,d}(t).
\end{equation}
Since Legendre polynomials are orthogonal with weight 
$(1-t^2)^{\frac{d-2}{2}}$ we have
\begin{equation}
\label{otn}
c_{jk}=\frac{\int_{-1}^1 (1-t)^j L_{k,d}(t) (1-t^2)^{\frac{d-2}{2}}dt}
{\int_{-1}^1 L_{k,d}(t)^2 (1-t^2)^{\frac{d-2}{2}}dt}.
\end{equation}
The denominator of (\ref{otn}) is equal to (see [M\"u]):
$$
\frac{\operatorname{vol}(S^d)}{\operatorname{vol}(S^{d-1}) \mu_{k,d}}=
\frac{\Gamma(\frac{d}{2})\sqrt{\pi}}{\Gamma(\frac{d+1}{2}) \mu_{k,d}},
$$ 
where the last equality follows from (\ref{vol2}).
The numerator of (\ref{otn}) is computed
using the Rodrigues rule ([M\"u])and the following integral (see [Er]): 
$$
\int_{-1}^1(1+t)^{\frac{d}{2}+k-1} (1-t)^{\frac{d}{2}+j-1}=
\frac{2^{k+j+d-1}\Gamma(\frac{d}{2}+k)\Gamma(\frac{d}{2}+j)}
{\Gamma(k+j+d)},
$$
Finally we get:
$$
c_{jk}=\frac{(-1)^k 2^j \Gamma(j+\frac{d}{2}) j!}{(j-k)!(j+k+d-1)!}
\frac{(4\pi)^{d/2} \mu_{k,d}}{\operatorname{vol}(S^d)}
$$
Let us  substitute this into (\ref{leg0}) and further on into (\ref{mm}). 
Note that $L_{k,d}(\cos 0)=L_{k,d}(1)=1$ for all $k$ (see [M\"u]).
Taking into account (\ref{leg}) and (\ref{vol1}) we obtain (\ref{sp}) after 
some easy combinatorial transformations.
This completes the proof of Theorem \ref{s} for $\omega \ge 3n$.
\qed

As we mentioned in section 1.1, it follows from the proof of 
Theorems \ref{o} and \ref{e} that in fact one can take $\omega \ge 2n$ 
(see section 4.3). 
\section{Proofs of the identities}
\subsection{Proof of Theorem \ref{S1}}
In this section we follow [Z].
We will prove a more general statement:
\begin{equation}
\label{S1g}
\sum_{j=0}^{\omega} \frac{1}{(\omega-j)! (j+n)! (2j+1)}\sum_{k=-j}^j
\frac{(-1)^k (x+k)^{2j+2n}}{(j-k)!(j+k)!} = 0,
\end{equation}
for $x\in {\Bbb R}$ and $\omega\ge 2n$.
If $x=0$ we get the original $S^1$-identity.
Note that we have symmetrized the summation limits in the inner sum ---
this is equivalent to multiplying the left-hand side by factor $2$.
Our aim is to make (\ref{S1g}) hypergeometric, i.e. to represent it 
as a function
\begin{equation}
\label{hyp}
_2F_1(a,b;c;z)=\sum_{m=0}^{\infty}
\frac{(a)_m (b)_m}{(c)_m}\frac{z^m}{m!},
\end{equation}
where $(t)_m=t(t+1)\cdots(t+m-1)$, $(t)_0=1$. 
Let $Ef(x)=f(x+1)$ be the shift operator. Then we can rewrite (\ref{S1g}) as
\begin{multline}
\label{z1}
\sum_{j=0}^\omega \frac{(-1)^j}{(\omega-j)! (j+n)! (2j+1)!}
\sum_{p=0}^{2j} (-1)^{p} \binom{2j}{p}E^{p-j}x^{2j+2n}=\\
\sum_{j=0}^{\omega}\frac{(-1)^j}{(\omega-j)! (j+n)! (2j+1)!} 
(E^{1/2}-E^{-1/2})^{2j}x^{2j+2n}.
\end{multline}
Using Taylor theorem $E=e^D$ where $D$ is the differentiation operator 
(see [GKP]) we have:
$$
(E^{1/2}-E^{-1/2})^{2j}=(e^{D/2}-e^{-D/2})^{2j}=P(D)^{2j} D^{2j},
$$
where 
\begin{equation}
\label{PD}
P(D)=\frac{2 \sinh D/2}{D}=\frac{e^{D/2}-e^{-D/2}}{D}=1+\frac{D^2}{24}+O(D^4).
\end{equation}
Substituting this into the sum and applying $D^{2j}$ to $x^{2j+2n}$ we get:
\begin{multline}
\label{ff}
\sum_{j=0}^{\omega} \frac{(-1)^j (2j+2n)! P(D)^{2j} x^{2n}}
{(\omega-j)! (j+n)! (2j+1)! (2n)!}=\\
\frac{1}{\omega! n!} \,{_2F_1}(n+1/2,-\omega;3/2;P(D)^2)x^{2n}=\\
\frac{1}{\omega! n!} \,{_2F_1}(1-n,\omega+3/2;3/2;P(D)^2)(1-P(D)^2)^{\omega-
n+1}
x^{2n}.
\end{multline}
The first equality is obtained by representing the sum as a hypergeometric
series and the second equality follows from the Euler transformation 
(see [GKP]):
\begin{equation}
\label{Euler}
_2F_1(a,b;c;z)=(1-z)^{c-a-b} {_2F_1}(c-a,c-b;c;z).
\end{equation}
Note that on both sides we have in fact polynomials in $D$ 
since $-\omega \le 0$ and $1-n \le 0$ and therefore both hypergeometric series
are finite (otherwise they would not be well defined).

On the other hand, due to (\ref{PD}) we have
$$
(1-P(D)^2)^{\omega-n+1}=O(D^{2\omega-2n+2}),
$$
and hence 
$$
(1-P(D)^2)^{\omega-n+1}x^{2n}=0
$$
for $\omega\ge 2n$.
This completes the proof of the $S^1$-identity. \qed

Note that for $\omega=2n-1$  the identity (\ref{S1}) does not hold (see [Z])
and hence $2n$ is ``sharp'' as was mentioned in section 1.1.
\subsection{Proof of Theorem \ref{S3}}
As in the previous section, we symmetrize the inner sumation indices and 
prove that
\begin{equation*}
\sum_{j=0}^\omega \frac{(-1)^n \Gamma(\omega+5/2)}{(\omega-j)! (j+n)! (2j+3)} 
\sum_{l=-j-1}^{j+1} (-1)^l \frac{l^2 (l^2-1)^{j+n}}{(j+l+1)! (j-l+1)!}=
-\frac{\sqrt{\pi}}{4 \cdot n!},
\end{equation*}
for $n\ge 1$, $\omega \ge 2n$.

We transform the inner sum:
$$
\frac{1}{(2j+2)!}\sum_{l=-j-1}^{j+1} \binom{2j+2}{j-l+1} (l^2-1)^{j+n} l^2=
$$
$$
\frac{(-1)^{j+1}}{(2j+2)!} \sum_{p=0}^{2j+2} (-1)^p (p-j-1)! \binom{2j+2}{p}
((p-j-1)^2-1)^{j+n}.
$$
Let us substitute this to the initial expression changing the summation index 
$j\to j+1$. Denote $\omega'=\omega+1$, $n'=n-1$.
We have
\begin{multline}
\label{11}
\sum_{j=0}^{\omega'} 
\frac{(-1)^{n'+1} \Gamma(\omega'+3/2) (-1)^j}{(\omega'-j)! (j+n')! (2j+1)!} 
\sum_{p=0}^{2j}(-1)^p \binom{2j}{p} (p-j)^2 ((p-j)^2-1)^{j+n'}.
\end{multline}
Let us open the last bracket. We get:
$$
\sum_{r=0}^{j+n'} (-1)^r \binom{n'+j}{r}\sum_{p=0}^{2j} (-1)^p 
\binom{2j}{p}(p-j)^{2j+2n'-2r+2}.
$$
Note that (see (1.13) in [Go])
\begin{equation}
\label{vychet1}
\sum_{p=0}^{2j} (-1)^p \binom{2j}{p}(p-j)^{s}=0
\end{equation}
for $s<2j$  and
\begin{equation}
\label{vychet2}
\sum_{p=0}^{2j} (-1)^p \binom{2j}{p}(p-j)^{2j}=(2j)!
\end{equation}
Therefore non-zero contribution comes only from $2j+2n'-2r+2 \ge 2j$, i.e.
$r\le n'+1$.
This implies that  (\ref{11}) can be rewritten as
\begin{multline}
\label{12}
(-1)^{n'+1} \Gamma(\omega'+3/2)
\sum_{r=0}^{n'+1} (-1)^r\binom{n'+1}{r}\cdot\\
\sum_{j=1}^{\omega'}
\frac{ (-1)^j}{(\omega'-j)! 
(j+t-1)! (2j+1)!}
\sum_{p=0}^{2j}(-1)^p \binom{2j}{p} (p-j)^{2j+2t}, 
\end{multline}
where $t=n'-r+1$. Consider the last two sums:
\begin{equation}
\label{zeil}
\sum_{j=1}^{\omega'}
\frac{(-1)^j}{(\omega'-j)! (j+t-1)! (2j+1)!}\sum_{p=0}^{2j}(-1)^p
\binom{2j}{p} (p-j)^{2j+2t}
\end{equation}
Let us show that (\ref{zeil}) vanishes for $\omega'\ge 2t+1$ which is always
the case since $\omega\ge 2n$ and $t\le n'+1=n$).
We use Lemma \ref{ff1} (see section 5.1) taking $s=1$ in (\ref{bb}). 
Applying the same arguments as in the proof of Theorem \ref{S1}
we get that (\ref{zeil}) vanishes for $r<n'+1$.  
Therefore the only non-zero contribution to (\ref{12}) comes from $r=n'+1$.
Taking this into account and substituting (\ref{vychet2})
into (\ref{12}) we finally obtain:
\begin{multline*}
\frac{\Gamma(\omega'+3/2)}{(n'+1)! \,a'!} \sum_{j=1}^{\omega'} 
\frac{(-1)^j j}
{2j+1}\binom{\omega'}{j}=
-\frac{\Gamma(\omega'+3/2) \sqrt{\pi}}{4 (n'+1)!\,
\Gamma(\omega'+3/2)}=-\frac{\sqrt{\pi}}{4\cdot n!}.
\end{multline*}
which completes the proof of Theorem \ref{S3}. \qed
\section{Proofs of Theorems \ref{o} and \ref{e}}
\subsection{Proof of Theorem \ref{o}}
Denote  $z=k+\alpha$. The inner sum  in (\ref{sp}) is equal to:
\begin{multline*}
\sum_{z=\alpha}^{j+\alpha}(-1)^{z+\alpha}\frac{2z\,(z+\alpha-1)!}{(z-\alpha)! 
(2\alpha)!}\binom{2j+2\alpha}{j+\alpha+z}
(z^2-\alpha^2)^{j+n}=\\
2\cdot \sum_{z=\alpha}^{j+\alpha} \frac{(-1)^{z+\alpha}}{(2\alpha)!}
\prod_{\beta=0}^{\alpha-1} (z^2-\beta^2) \,
\binom{2j+2\alpha}{j+\alpha+z}(z^2-\alpha^2)^{j+n} =\\
\sum_{z=-j-\alpha}^{j+\alpha} \frac{(-1)^{z+\alpha}}{(2\alpha)!} 
\prod_{\beta=0}^{\alpha-1}
(z^2-\beta^2) \binom{2j+2\alpha}{j+\alpha+z} (z^2-\alpha^2)^{j+n}.
\end{multline*}
Denote $l=j+\alpha+z$. Then the last sum can be rewritten as
\begin{equation}
\label{pr}
\frac{(-1)^j}{(2\alpha)!}\sum_{l=0}^{2j+2\alpha} (-1)^l 
\prod_{\beta=0}^{\alpha-1} ((l-j-\alpha)^2-\beta^2)\, 
\binom{2j+2\alpha}{l} ((l-j-\alpha)^2-\alpha^2)^{j+n}
\end{equation}
Let $\omega'=\omega+\alpha$, $n'=n-\alpha$ and let $j:=j+\alpha$ be the new
summation index. Due to (\ref{pr})  we can represent (\ref{sp}) as:
\begin{multline}
\label{!!}
a_{n,2\alpha+1}=
\frac{2(-1)^{n'}\Gamma(\omega'+\frac{3}{2})}{(2\alpha)!}
\sum_{s=1}^{\alpha}K_s^{\alpha} 
\cdot\\
\sum_{j=0}^{\omega'}\frac{(-1)^{j}}{(\omega'-j)!(j+n'-r)!(2j+1)!}\cdot\\
\sum_{r=0}^{j+n'}\frac{(-1)^r\alpha^{2r}}{r!}
\sum_{l=0}^{2j}(-1)^l \binom{2j}{l}(l-j)^{2j+2n'-2r+2s},
\end{multline}
where $K_s^\alpha$ are defined by (\ref{cno}).
Note that if $2j+2n'-2r+2s<2j$ the last sum vanishes due to (\ref{vychet1}).
Therefore if $r\le n'+s$ we can rewrite (\ref{!!}) as 
\begin{multline}
\label{!!1}
a_{n,2\alpha+1}=
\frac{2(-1)^{n'}\Gamma(\omega'+\frac{3}{2})}{(2\alpha)!}
\sum_{s=1}^{\alpha}K_s^{\alpha} 
\sum_{r=0}^{n'+s}\frac{(-1)^r\alpha^{2r}}{r!}
\cdot\\
\sum_{j=0}^{\omega'}\frac{(-1)^{j}}{(\omega'-j)!(j+n'-r)!(2j+1)!}\cdot\\
\sum_{l=0}^{2j}(-1)^l \binom{2j}{l}(l-j)^{2j+2n'-2r+2s},
\end{multline}
using the fact that $(j+n'-r)!=0$ for $r>j+n'$.
Let us note that  Lemma \ref{ff1} implies that the last two sums 
in (\ref{!!1}) vanish if $r<n'+s$ and $\omega \ge 2n$. Indeed, this follows 
from  (\ref{bb})  for $t=n'-r+s$ in the same way as vanishing of 
(\ref{zeil}) in the proof of Theorem \ref{S3}. Therefore the only non-zero 
contribution again comes only from $r=n'+s$ when the inner sum is equal 
to $(2j)!$ by (\ref{vychet2}). 
Hence we obtain:
\begin{multline}
\label{grr}
a_{n,2\alpha+1}=\frac{2\Gamma(\omega'+\frac{3}{2})}{(2\alpha)!}
\sum_{s=1}^{\alpha} K_s^{\alpha} \frac{(-1)^{s}\alpha^{2n'+2s}}{(n'+s)!}
\sum_{j=0}^{\omega'}\frac{(-1)^j}{(\omega'-j)!(j-s)!(2j+1)}.
\end{multline}
Note that 
\begin{multline*}
\sum_{j=0}^{\omega'}\frac{(-1)^j}{(\omega'-j)!(j-s)!(2j+1)}=\\
\frac{(-1)^s}{(\omega'-s)!}\sum_{j=0}^{\omega'-s} (-1)^j \binom{\omega'-s}{j}
\frac{1}{2j+2s+1}=\\
\frac{(-1)^s}{(\omega'-s)!}\int_0^1 
\left(\sum_{j=0}^{\omega'-s} (-1)^j \binom{\omega'-s}{j}x^{2j+2s}\right)dx=\\
\frac{(-1)^s}{(\omega'-s)!}\int_0^1 x^{2s} (1-x^2)^{\omega'-s}dx=
\frac{(-1)^s \Gamma(s+\frac{1}{2})}{2 \Gamma(\omega'+3/2)},
\end{multline*}
where the last equality follows from ([GR]).
Substituting this into (\ref{grr}) after certain cancellations we 
obtain (\ref{odd}).
In particular, taking $\alpha=2$ and $\alpha=3$ we get (\ref{57}).
The proof of Theorem \ref{o} is complete. \qed
\subsection{Proof of Theorem \ref{e}} The first steps of the proof
are similar to that of Theorem \ref{o}.
Let $n'=n-\nu+1$, $\omega'=\omega+\nu-1$ and let  
$j:=j+\nu-1$ be the new summation index.
Similarly to (\ref{!!1}) we obtain the following formula from (\ref{sp}):
\begin{multline}
\label{!!!}
a_{n,2\nu}=
\frac{2(-1)^{n'+1}(\omega'+1)!}{(2\nu-1))!}
\sum_{s=0}^{\nu-1}K_s^{\nu} \sum_{r=0}^{n'+s}
\frac{(-1)^r(\nu-\frac{1}{2})^{2r}}{r!}\cdot\\
\sum_{j=0}^{\omega'}\frac{(-1)^{j}}{(\omega'-j)!(j+n'-r)!(2j+2)!}\cdot\\
\sum_{l=0}^{2j+1}(-1)^l \binom{2j+1}{l}(l-j-\frac{1}{2})^{2j+2n'-2r+2s+1},
\end{multline}
However, from this moment the situation is quite different. 
If in the proof of Theorem \ref{o} only one term  corresponding to 
$r=n'+s$ gave a non-zero contribution, now all terms with $n'-r+s \ge 0$
contribute to the sum. 
Indeed, repeating the arguments of Theorem \ref{S3}
we get that the last two sums in (\ref{!!!}) 
are equal to: \begin{equation}
\label{haha}
\sum_{j=0}^{\omega'} \frac{(-1)^j (2j+2n'-2r+2s+1)! P^{2j+1}}
{(\omega'-j)! (j+n'-r)! (2j+2)! (2n'-2r+2s)!} x^{2n'-2r+2s}|_{x=0}
\end{equation}
where $P$ is given by (\ref{PD}).
Setting $t=n'-r+s$ in (\ref{e2}) in Lemma \ref{ff2} (see section 5.2)
we get that if $\omega \ge 2n$,  (\ref{haha}) is equal to
\begin{equation}
\label{bern0}
\frac{(n'-r)_s}{2 (n'+s-r)! (\omega'+1)!}
P^{-1}(x^{2n'-2r+2s})|_{x=0}
\end{equation}
Let us compute $P^{-1}(x^{2t})|_{x=0}$.
We have 
$$P^{-1}=\frac{D}{e^{-D/2}-e^{D/2}}=\sum_{i=0}^{\infty}P_{2i}D^{2i}$$
and
$$P^{-1}(x^{2t})|_{x=0}=(2t)! P_{2t}.$$
Computing $P_{2t}$ we get Bernoulli numbers. Indeed,
\begin{equation}
\label{bern1}
(2t)! P_{2t}=-2\left(\frac{B_{2t}}{2^{2t}}-\frac{B_{2t}}{2}\right)
\end{equation}
Indeed, by a well-known formula (see [GKP])
$$
\frac{z}{e^{z}-1}=\sum_{n=0}^\infty \frac{B_n z^n}{n!},
$$
and on the other hand
$$
\frac{z/2}{e^{z/2}-1}-\frac{1}{2}\frac{z}{e^{z}-1}=-\frac{1}{2}P^{-1}(z)
$$
which implies (\ref{bern1}).
Let us substitute (\ref{bern1}) into (\ref{bern0}) and further into 
(\ref{!!!}). After certain combinatorial transformations
we obtain (\ref{even}). In particular, we take $t$ as the new summation
index, $t=0,1,\dots, n-\nu+1+s$. Note that if $n<\nu$ then $t\le s$
and hence $(t-s)_s=(t-s)(t-s+1)\cdots (t-1)=0$ unless $t=0$ when 
$(-s)_{s}=(-1)^s s!$. This explains why the second sum disappears in 
(\ref{even}) for $n<\nu$.

It is easy to check that taking $\nu=1$ we get (\ref{S2}). 
The proof is complete.
\qed
\subsection{Proof of Theorem \ref{s} for $2n\le\omega<3n$}
Take the arguments in the proofs of Theorems \ref{o} and \ref{e}
in the reverse order. Starting with (\ref{odd}) in odd dimensions and 
(\ref{even}) in even dimensions we arrive to (\ref{sp}).
Note that the proofs of Theorems \ref{o} and \ref{e} are valid for
$\omega \ge 2n$ (cf. Theorems \ref{S1} and \ref{S3}) 
and hence formula (\ref{sp}) holds under the same condition.
This completes the proof of Theorem \ref{s}. \qed
\section{Auxiliary combinatorial lemmas}
\subsection{Odd dimensions}
\begin{lemma} Let $\omega'\ge 2t+s$, $s\ge 0$, $t\ge 1$. Then
\label{ff1}
\begin{multline}
\label{aa}
\sum_{j=0}^{\omega'} \frac{(-1)^j(2j+2t)! z^j}{(\omega'-j)!(j+t-s)!(2j+1)!}=\\
\sum_{k=0}^s Q_{k,s}(z)\, _2F_1(-\omega'+k,\frac{1}{2}+t+k;\frac{3}{2}+k;z),
\end{multline}
where $Q_{k,s}(z)$ are some polynomials in $z$.
Moreover,
\begin{equation}
\label{bb}
\sum_{j=0}^{\omega'} \frac{(-1)^j (2j+2t)!}{(\omega-j)! (j+t-s)! (2j+1)!}
P(D)^{2j} x^{2t}=0,
\end{equation}
where $P(D)$ is defined by (\ref{PD}).
\end{lemma}
\noindent{\bf Proof.} 
Denote the sum at the left hand side by $\sigma_s(z)$.
Let us proceed by induction. For $s=0$ the statement follows from (\ref{ff}).
Suppose we proved it for all $s\le s_0$. Let us prove it for $s_0+1$.
It is easy to see that
\begin{equation}
\label{sig}
\sigma_{s_0+1}(z)=(t-s_0)\sigma_{s_0}(z)+z \frac{d \sigma_{s_0}}{dz}
\end{equation}
By the induction hypothesis and the rule for differentiation of a 
hypergeometric function (see [Er]) we obtain:
\begin{multline*}
\frac{d\sigma_{\sigma_0}}{dz}=\sum_{k=0}^{s_0} Q_{k,s_0}'(z)
\, _2F_1(-\omega'+k,\frac{1}{2}+t+k;\frac{3}{2}+k;z)+\\
\sum_{k=0}^{s_0} 
Q_{k,s_0}(z) 
\frac{(k-\omega')(\frac{1}{2}+t+k)}{k+\frac{3}{2}} \, _2F_1(-\omega'+k+1,
\frac{3}{2}+t+k;\frac{5}{2}+k;z)
\end{multline*}
Substituting this to (\ref{sig}) implies (\ref{aa}).

Let us prove (\ref{bb}). We use (\ref{aa}) and  apply arguments of the 
previous section starting with (\ref{ff}) to each term of the sum 
$\sigma_s(z)$. Note that each hypergeometric function in the right-hand
side of (\ref{aa}) is in fact a finite series since $-\omega'+k <0$ for
all $k=0,1,\dots,s$.
Due to (\ref{Euler}) we have the 
following condition for vanishing of the left hand side in (\ref{bb}):
\begin{equation}
\label{imp}
3/2+k-1/2-t-k+\omega'-k=\omega'-k-t+1>t,
\end{equation}
that is $\omega'\ge 2t+k$. But we have supposed that $\omega'\ge 2t+k$
and since $k\le s$ we get (\ref{imp}). The last thing we have to
verify is that using the Euler transformation (\ref{Euler}) we always get
a finite hypergeometric series. This is indeed so since 
$$3/2+k-1/2-t-k=1-t\le 0$$ due to the condition $t\ge 1$.
This completes the proof of the Lemma. \qed
\subsection{Even dimensions}
\begin{lemma}
\label{ff2}
Let $\omega' \ge 2t+s$, $s\ge 0$, $t\ge 0$. Then 
\begin{multline}
\label{e1}
\sum_{j=0}^{\omega'} \frac{(-1)^j (2j+2t+1)! z^{j+1}}{(\omega'-j)! (j+t-s)! 
(2j+2)!}=\\
\frac{(2t)! (t-s)_s}{2 (\omega'+1)! \, t!}+
\sum_{k=0}^s Q_{k,s}(z){_2F_1}(-1-\omega'+k, \frac{1}{2}+t+k; \frac{1}{2}+k;z),
\end{multline}
where $Q_{k,s}(z)$ are some polynomials in $z$.
Moreover,
\begin{multline}
\label{e2}
\sum_{j=0}^{\omega'} \frac{(-1)^j (2j+2t+1)! P^{2j+1}}
{(\omega'-j)! (j+t-s)! (2j+2)!} x^{2t}|_{x=0}
=\frac{(2t)! (t-s)_s}{2 (\omega'+1)! \, t!}
P^{-1}(x^{2t})|_{x=0},
\end{multline}
where $P(D)$ is defined by (\ref{PD}).
\end{lemma}
\noindent{\bf Proof.} 
Again, we proceed by induction over $s$. 
For $s=0$ this can be checked by a direct computation (e.g. using [W]). 
Denoting the left-hand side of (\ref{e1}) by $\zeta_s(z)$
similarly to (\ref{sig}) we have
\begin{equation}
\label{sig1}
\zeta_{s_0+1}(z)=(t-s_0-1)\zeta_{s_0}(z)+z \frac{d \zeta_{s_0}}{dz}
\end{equation}
As in the proof of Lemma \ref{ff1} this implies the induction step and
proves (\ref{e1}). The relation  (\ref{e2}) follows from  (\ref{e1})
in a similar way as (\ref{bb}) follows from (\ref{aa}).
\qed
\section*{References}

\noindent [Be] M. Berger, Geometry of the spectrum, Proc. Symp. Pure Math. 27
(1975), 129-152.

\smallskip

\noindent [CW] R.S. Cahn, J.A. Wolf, Zeta functions and their asymptotic
expansions for compact symmetric spaces of rank one, Comment. Math. Helv.,
51 (1976),1-21.

\smallskip

\noindent [Ca] R. Camporesi, Harmonic analysis and propagators on 
homogeneous spaces, Phys. Rep. 196 (1990), no. 1-2, 1-134.

\smallskip

\noindent [ELV] E. Elizalde, M. Lygren, D.V. Vassilevich, Antisymmetric
tensor fields on spheres: functional determinants and non-local counterterms,
J. Math. Phys. 37 (1996), no. 7, 3105-3117.

\smallskip

\noindent [Er] A. Erd\'elyi et. al., Higher transcendental functions, vol. 1,
McGraw-Hill, 1953.

\smallskip

\noindent [DK] J.S. Dowker, K. Kirsten, Spinors and forms on the ball and the
generalized cone, Comm. Anal. Geom. 7 (1999), no. 3, 641-679.

\smallskip

\noindent [Gi] P. Gilkey, The index theorem and the heat equation, Math. 
Lect.Series,  Publish or Perish, 1974.

\smallskip

\noindent [Go]  H.W. Gould, Combinatorial identities, 1972.

\smallskip

\noindent [GR] I.S. Gradshtein, I.M. Ryzhik, Table of integrals, series and
products, Academic Press,  (1980).

\smallskip

\noindent [GKP] R. Graham, D. Knuth, O. Patashnik, Concrete Mathematics.
A foundation for computer Science, Addison-Wesley, 1994.

\smallskip

\noindent [MS]  H.P. McKean, Jr., I.M. Singer,
Curvature and the eigenvalues of the
Laplacian, J. Diff. Geom. 1 (1967), 43-69.

\smallskip

\noindent [Mi] J. Milnor, Morse Theory, Princeton University Press, 1963.

\smallskip

\noindent [M\"u] C. M\"uller, Analysis of spherical symmetries in 
Euclidean spaces, Applied Mathematical Sciences, 129, Springer-Verlag, 1998.

\smallskip

\noindent [P] I. Polterovich, Heat invariants of Riemannian manifolds,
Israel J. Math. 119 (2000), 239-252.

\smallskip

\noindent [Wo] S. Wolfram, Mathematica: a system for doing mathematics
by computer, Addison -- Wesley, 1991.

\smallskip

\noindent [We] G. Weingart, A characterization of the heat kernel coefficients,
math.DG/0105144.

\smallskip

\noindent [Z]. D. Zeilberger, Proof of an identity conjectured by Iossif
Polterovitch that came up in the Agmon-Kannai asymptotic theory of the
heat kernel, (2000), http://www.math.temple.edu/~zeilberg/pj.html.
\end{document}